# A Product to Sum Approach for Matrix Filling in a Hierarchical Finite-Element Method

Ehsan Khodapanah

*Abstract*−**This paper presents a product to sum approach for a fast and efficient matrix filling in a hierarchical finite-element method (FEM). Due to the existence of a coupling factor arising from the material and Jacobian inhomogeneities in curved inhomogeneous elements, the calculation of the FEM matrix elements should be carried out through full multidimensional integrations. This reduces the efficiency of the higher order FEM solvers especially when the coupling factor varies rapidly inside the elements. In the product to sum approach, every product of the basis and weighting polynomials is replaced with a sum of appropriate polynomials. This reduces the number of required multidimensional integrals significantly and converts the integration into a summation. Therefore, the method will be efficient if the number of summation terms in the product to sum conversion is as small as possible.**

*Index terms*− **Product to sum, finite-element method, matrix filling.**

## I.  INTRODUCTION

The finite-element method (FEM) is a robust and accurate numerical method for solving vector partial differential equations (PDEs) appearing in a variety of research areas especially in the area of electromagnetic field computations [1]-[3].

When the FEM is applied on curved elements containing inhomogeneous materials, the FEM matrix elements should be calculated through full multidimensional integrations, which are costly and time-consuming. The coupling factor in the matrix element integrals is due to the material inhomogeneity or Jacobian inhomogeneity or both. In the case of straight triangular and tetrahedral elements with homogeneous materials, the coupling factor is reduced to a constant and is removed from the integral sign. Therefore, a universal matrix approach can be applied for a fast matrix filling [4]. However, in general even in a straight quadrilateral or hexahedral element with a homogeneous material, the coupling factor is not a constant and hence the universal matrix approach is not applicable.

The generalized universal matrix approach [5] and a universal array approach [6] have also been proposed and can be used for fast matrix filling in a curved inhomogeneous element. These methods are based on a polynomial approximation of the coupling factor and conversion of the integral into a sum over universal integrals that can be calculated and stored beforehand. However, these methods are more efficient than direct numerical integrations if the polynomial order of the field expansion inside the element is chosen larger than the half of the order of the coupling factor approximation. The efficiency of these methods increases as the difference between the polynomial orders of the field and coupling factor is increased.

Separable Jacobian methods have also been proposed for fast matrix filling in the 2-D curved homogeneous elements [7]-[9]. These methods are based on the fact that in the ring-type elements the coupling factor is separable [7], [8]. On the other hand, an arbitrary curved polygonal element can be divided into special curved triangles for which the coupling factor is separable [9]. Therefore, 2-D integrals of the FEM matrix elements are reduced to 1-D ones that can be calculated very quickly.

In this paper, we present a product to sum approach for fast matrix filling in a hierarchical curved inhomogeneous quadrilateral or hexahedral element. The method is based on the fact that the product of two polynomials can be represented as a sum of appropriate polynomials. Therefore, the multidimensional integrals are reduced to summations over another type of integrals. In this procedure, the number of multidimensional integrals is reduced significantly. However, the efficiency of the method depends on the number of terms in the equivalent polynomial representations of the products.

The author is with the Sum Institute for Computational Physics, Tabriz, Iran (e-mail: ehskhodapanah@ yahoo.com).

## II. General representation of the product to sum approach

We consider a second order vector PDE in a 3-D space. The unknown vector field, governed by the vector PDE, is denoted by **A**. In a hierarchical version of the FEM, different components of the vector field **A** in a curved hexahedral element can be expanded in terms of the following polynomial basis functions

$$A_\gamma^{mnp} = P_m^\gamma(u)Q_n^\gamma(v)R_p^\gamma(w) \quad \begin{cases} m = 0,1,\cdots N_u^\gamma - 1 \\ n = 0,1,\cdots N_v^\gamma - 1 \\ p = 0,1,\cdots N_w^\gamma - 1 \\ \gamma = \gamma_1, \gamma_2, \gamma_3 \end{cases} \quad (1)$$

where $\gamma_1$, $\gamma_2$, and $\gamma_3$ stand for the scalar components of the vector field, $\{P_m^\gamma(u)\}_{m=0}^{N_u^\gamma-1}$, $\{Q_n^\gamma(v)\}_{n=0}^{N_v^\gamma-1}$, and $\{R_p^\gamma(w)\}_{p=0}^{N_w^\gamma-1}$ are three sets of appropriate orthogonal polynomials in one dimension and $(u,v,w)$ are the coordinates of the reference element. By applying a standard Galerkin method to discretize the vector PDE, we arrive at the following general expressions for the matrix elements

$$M_{\gamma_i\gamma_j}^{m_1m_2n_1n_2p_1p_2} = \sum_{s=1}^{N_s} \iiint\limits_{-1}^{1} \begin{bmatrix} P_{m_1}^{\gamma_i}{'} \\ P_{m_1}^{\gamma_i} \end{bmatrix} \begin{bmatrix} P_{m_2}^{\gamma_j}{'} \\ P_{m_2}^{\gamma_j} \end{bmatrix} \begin{bmatrix} Q_{n_1}^{\gamma_i}{'} \\ Q_{n_1}^{\gamma_i} \end{bmatrix} \begin{bmatrix} Q_{n_2}^{\gamma_j}{'} \\ Q_{n_2}^{\gamma_j} \end{bmatrix} \begin{bmatrix} R_{p_1}^{\gamma_i}{'} \\ R_{p_1}^{\gamma_i} \end{bmatrix} \begin{bmatrix} R_{p_2}^{\gamma_j}{'} \\ R_{p_2}^{\gamma_j} \end{bmatrix} \alpha_s^{\gamma_i\gamma_j}(u,v,w)dudvdw \quad (2)$$

where $i,j = 1,2,3$, $\alpha_s^{\gamma_i\gamma_j}$ is the coupling factor, originating from the material inhomogeneity and the Jacobian components, and one element in every column vector in (2) should be chosen depending upon $s$, $\gamma_i$, and $\gamma_j$. In general $\alpha_s^{\gamma_i\gamma_j}$ is not separable and hence a full 3-D integration should be performed to obtain every 3-D integral inside the sum in the right-hand-side (RHS) of (2). This means that for every submatrix in the RHS of (2) that has a contribution to $\left[M_{\gamma_i\gamma_j}\right]$, $N_u^{\gamma_i}N_u^{\gamma_j}N_v^{\gamma_i}N_v^{\gamma_j}N_w^{\gamma_i}N_w^{\gamma_j}$ full 3-D integrals should be evaluated. In the product, to sum approach we replace every product of two polynomials of the same variable with a sum over a set of appropriate polynomials. For example

$$\begin{bmatrix} P_{m_1}^{\gamma_i}{'}(u) \\ P_{m_1}^{\gamma_i}(u) \end{bmatrix}_s \begin{bmatrix} P_{m_2}^{\gamma_j}{'}(u) \\ P_{m_2}^{\gamma_j}(u) \end{bmatrix}_s = \sum_{k_1=0}^{K_1} a_{k_1}^{s\gamma_i\gamma_j} P_{k_1}^{s\gamma_i\gamma_j}(u) \quad (3)$$

where $K_1 = m_1 + m_2$ or $m_1 + m_2 - 1$ or $m_1 + m_2 - 2$ depending on whether the polynomial or its derivative is considered and $P_{k_1}^{\gamma_1}$ may be chosen one among $P_{k_1}^{\gamma_1}$, $P_{k_1}^{\gamma_2}$, and $P_{k_1}^{\gamma_3}$ or even another polynomial if appropriate.

By substituting the product to sum formulas (3) into (2), we obtain the following expressions for the elements of $\left[M_{\gamma_i\gamma_j}\right]$

$$M_{\gamma_i\gamma_j}^{m_1m_2n_1n_2p_1p_2} = \sum_{s=1}^{N_s} \sum_{k_1,k_2,k_3} a_{k_1k_2k_3}^{s\gamma_i\gamma_j} \iiint\limits_{-1}^{1} P_{k_1}^{s\gamma_i\gamma_j}(u)Q_{k_2}^{s\gamma_i\gamma_j}(v)R_{k_3}^{s\gamma_i\gamma_j}(w)\alpha_s^{\gamma_i\gamma_j}(u,v,w)dudvdw \quad (4)$$

Now it is sufficient to calculate only $\left(N_u^{\gamma_i} + N_u^{\gamma_j}\right)\left(N_v^{\gamma_i} + N_v^{\gamma_j}\right)\left(N_w^{\gamma_i} + N_w^{\gamma_j}\right)$ 3-D integrals to construct every submatrix that contributes to $\left[M_{\gamma_i\gamma_j}\right]$. However, the product to sum approach is efficient if total number of terms in the triple sum in (4) is relatively small compared to the total number of sample points in the direct 3-D numerical integration of (4). In other words, the efficiency of the method is directly proportional to the sparsity of the vector containing the coefficients of the polynomials in the RHS of the product to sum formulas in (3). Therefore, the key question is that are there any polynomials that lead to a sparse RHS in the product to sum formulas? The answer is yes; one such choice is the Chebyshev polynomials that are used in the next section.

## III. A SPECIFIC EXAMPLE

In order to clarify the method outlined in section II and for the ease of computer implementation, we consider a two-dimensional curl-curl equation in the area of electromagnetic field computation as a specific example. Indeed we solve the following vector PDE

$$\nabla \times \frac{1}{\mu_r} \nabla \times \mathbf{E} - k_0^2 \varepsilon_r \mathbf{E} = 0 \qquad (5)$$

The domain of the problem is a fourth-order curved quadrilateral domain partially filled with a continuously-varying inhomogeneous material (Fig. 1). As in the standard FEM, we first divide the domain of Fig. 1 into a number of smaller curved quadrilateral elements. Inside each element, the vector field components are expanded in terms of the following basis functions

$$\begin{cases} E_u = U_m(u)\, T_n(v) & \begin{cases} m = 0,1,\cdots,M-1 \\ n = 0,1,\cdots,N \end{cases} \\ E_v = T_m(u)\, U_n(v) & \begin{cases} m = 0,1,\cdots,M \\ n = 0,1,\cdots,N-1 \end{cases} \end{cases} \qquad (6)$$

After applying a standard Galerkin method to (5), we obtain the following general expressions for the stiffness and mass matrices

$$S_{tb} = \iint \frac{1}{\mu_r} \nabla \times \mathbf{E}_t \cdot \nabla \times \mathbf{E}_b \; dxdy$$
$$M_{tb} = \iint \varepsilon_r \mathbf{E}_t \cdot \mathbf{E}_b \; dxdy \qquad (7)$$

where $t$ and $b$ stand for the testing and basis functions, respectively, and $[S]$ and $[M]$ are the stiffness and mass matrices, respectively. Substitution of (6) into (7) leads to the following explicit expressions for the stiffness and mass submatrices

$$S_{uu}^{m_1 m_2 n_1 n_2} = n_1 n_2 \iint_{-1}^{1} U_{m_1}(u) U_{m_2}(u) U_{n_1-1}(v) U_{n_2-1}(v) \frac{dudv}{\mu_r J} \qquad (8)$$

$$S_{uv}^{m_1 m_2 n_1 n_2} = -n_1 m_2 \iint_{-1}^{1} U_{m_1}(u) U_{m_2-1}(u) U_{n_1-1}(v) U_{n_2}(v) \frac{dudv}{\mu_r J} \qquad (9)$$

$$S_{vu}^{m_1 m_2 n_1 n_2} = -m_1 n_2 \iint_{-1}^{1} U_{m_1-1}(u) U_{m_2}(u) U_{n_1}(v) U_{n_2-1}(v) \frac{dudv}{\mu_r J} \qquad (10)$$

$$S_{vv}^{m_1 m_2 n_1 n_2} = m_1 m_2 \iint_{-1}^{1} U_{m_1-1}(u) U_{m_2-1}(u) U_{n_1}(v) U_{n_2}(v) \frac{dudv}{\mu_r J} \qquad (11)$$

$$M_{uu}^{m_1 m_2 n_1 n_2} = \iint_{-1}^{1} U_{m_1}(u) U_{m_2}(u) T_{n_1}(v) T_{n_2}(v) \frac{\varepsilon_r (x_v^2 + y_v^2)}{J} dudv \qquad (12)$$

$$M_{uv}^{m_1 m_2 n_1 n_2} = -\iint_{-1}^{1} U_{m_1}(u) T_{m_2}(u) T_{n_1}(v) U_{n_2}(v) \frac{\varepsilon_r (x_u x_v + y_u y_v)}{J} dudv \qquad (13)$$

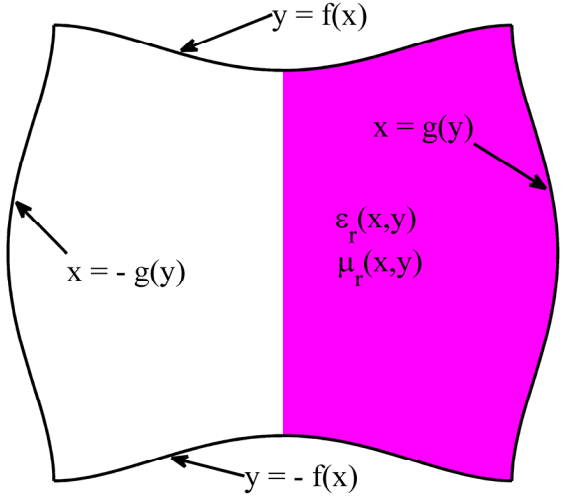

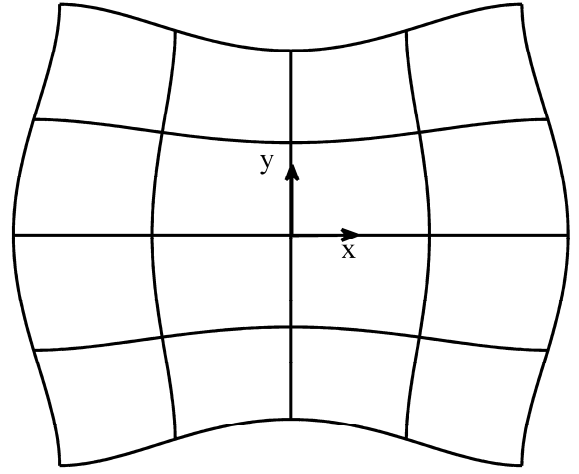

Fig. 2. The finite-element mesh for the domain of Fig. 1.

Fig. 1. Geometry of the 2-D problem considered in section III.

$f(x) = -0.2(x^2-1)^2 + 1,\ g(y) = 0.2(y^2-1)^2 + 1,\ \mu_r = 1,$ and
$\varepsilon_r = 2\exp(x+y+2).$

$$M_{vu}^{m_1 m_2 n_1 n_2} = -\iint\limits_{-1}^{1} T_{m_1}(u)U_{m_2}(u)U_{n_1}(v)T_{n_2}(v)\frac{\varepsilon_r(x_u x_v + y_u y_v)}{J}\,du\,dv \quad (14)$$

$$M_{vv}^{m_1 m_2 n_1 n_2} = \iint\limits_{-1}^{1} T_{m_1}(u)T_{m_2}(u)U_{n_1}(v)U_{n_2}(v)\frac{\varepsilon_r(x_u^2 + y_u^2)}{J}\,du\,dv \quad (15)$$

As can be realized from (8)-(15) every submatrix element is obtained through a 2-D integration in the reference element. As mentioned in the foregoing section, in the product to sum approach, we replace every product term containing the product of two chebyshev polynomials of a same variable with a sum of two appropriate polynomials by applying the following formulas

$$\begin{cases} U_{m_1}(u)U_{m_2}(u) = T_{|m_1-m_2|}^{ns}(u) - T_{m_1+m_2+2}^{ns}(u) \\[2mm] T_{m_1}(u)T_{m_2}(u) = T_{|m_1-m_2|}(u) + T_{m_1+m_2}(u) \\[2mm] U_{m_1}(u)T_{m_2}(u) = U_{m_1-m_2}(u) + U_{m_1+m_2}(u) \end{cases} \quad (16)$$

where the nonsingular Chebyshev polynomials are defined as

$$T_n^{ns}(u) = \begin{cases} \dfrac{T_n(u)-1}{2(1-u^2)} & : \ n \text{ is even} \\[3mm] \dfrac{T_n(u)-u}{2(1-u^2)} & : \ n \text{ is odd} \end{cases} \quad (17)$$

For example, for the $[M_{uu}]$ elements we have

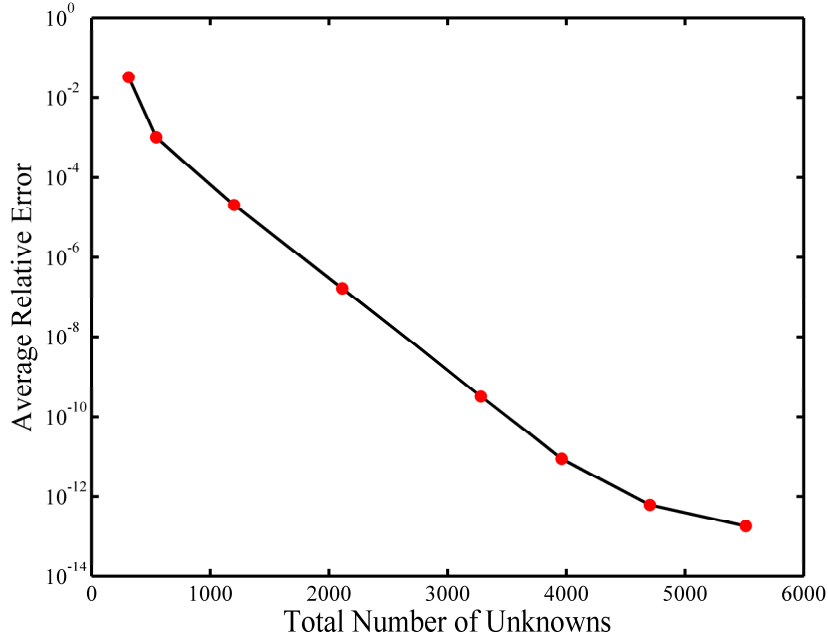

Fig. 3. Convergence pattern for the problem of Fig. 1.

Table I computational times for the product to sum approach and the direct numerical integration approach in solving the problem of Fig. 1.

| | Product to sum approach | | Direct numerical integration | | Reduction factor in the matrix filling time | |
|---|---|---|---|---|---|---|
| | Matrix filling time (s) | Total computational time (s) | Matrix filling time (s) | Total computational time (s) | From the left columns | Predicted value (section IV) |
| $M = N = 3$ | 0.104 | 0.33 | 0.100 | 0.33 | 0.96 | 0.6 |
| $M = N = 4$ | 0.121 | 0.38 | 0.148 | 0.41 | 1.22 | 1 |
| $M = N = 6$ | 0.170 | 0.5 | 0.428 | 0.76 | 2.5 | 2.4 |
| $M = N = 8$ | 0.229 | 0.8 | 1.103 | 1.66 | 4.8 | 4.25 |
| $M = N = 10$ | 0.328 | 1.32 | 2.375 | 3.35 | 7.25 | 6.7 |
| $M = N = 12$ | 0.461 | 2.2 | 4.65 | 6.36 | 10 | 9.6 |
| $M = N = 14$ | 0.634 | 3.6 | 8.21 | 11.2 | 13 | 13 |
| $M = N = 16$ | 0.934 | 5.9 | 14 | 19 | 15 | 17 |
| $M = N = 18$ | 1.29 | 8.8 | 22.14 | 29.7 | 17 | 21.6 |

$$M_{uu}^{m_1 m_2 n_1 n_2} = \iint\limits_{-1}^{1} \left[ T_{|m_1 - m_2|}^{ns}(u) T_{|n_1 - n_2|}(v) + T_{|m_1 - m_2|}^{ns}(u) T_{n_1 + n_2}(v) - T_{m_1 + m_2 + 2}^{ns}(u) T_{|n_1 - n_2|}(v) \right.$$
$$\left. - T_{m_1 + m_2 + 2}^{ns}(u) T_{n_1 + n_2}(v) \right] \frac{\varepsilon_r (x_v^2 + y_v^2)}{2J} \, du \, dv \qquad (18)$$

(18) shows that a four-term sum can be used to calculate the $[M_{uu}]$ elements instead of a direct numerical integration containing a large number of sums and products. However, the simple formula in (18) requires calculate the following 2-D integrals as a starting point

$$\iint\limits_{-1}^{1} T_m^{ns}(u)T_n(v)\frac{\varepsilon_r(x_v^2+y_v^2)}{2J}\,du\,dv \qquad (19)$$

Total number of 2-D integrals in (19) is approximately $4D$ where $D \approx MN$, while total number of 2-D integrals in the direct numerical calculation of $[M_{uu}]$ is approximately $D^2/4$ (taking the four-fold symmetry of $[M_{uu}]$ into account). Hence, the reduction factor for the number of 2-D integrals is $D/16$.

## IV. OPERATION COUNT ESTIMATION FOR THE PROPOSED METHOD

This section, we would like to give a relatively more accurate number for the arithmetic operation count of the product to sum approach applied to the problem of section III. Before proceeding, we make some assumptions to simplify the analysis. We assume that the curved elements in the domain decomposition are isoparametric and the geometrical mapping into the reference element is Lagrangian. These assumptions allow to find an upper limit for the number of operations required for the calculation of the Jacobian components. We also assume that the total number of sampling points for a direct 2-D integration is $2D$, which is a lower limit for an accurate numerical integration when the Jacobian and material inhomogeneities vary rapidly.

Total number of multiply-add operations for the calculation of the $x$, $y$, and four Jacobian components is $12D^2$. Total number of multiply-add operations for the calculation of the 2-D integrals as starting points is $32D^2$. Total number of add operations due to the product to sum formulas (e.g., (18)) is $3 \times 1.75D^2$. Total number of multiply operations for the calculation of the stiffness submatrices is $2 \times 1.5D^2$. Finally, there is a need for a rearrangement of the Chebyshev polynomials of the first kind in (6) to ensure the curl-conformity of the basis functions. This final step requires $9D^2$ add operations. Therefore, the total number of multiply-add operations in the product to sum formula is approximately $53D^2$. On the other hand, the total number of multiply-add operations in the direct numerical calculations is $2D \times 1.75D^2 = 3.5D^3$ (neglecting the Jacobian calculations and the multiply operations in the stiffness submatrices and the rearrangements).

The above analysis shows that the number of arithmetic operations is reduced by a factor of $D/15$ in the proposed method compared with the direct numerical integration.

## V. NUMERICAL RESULTS

In this section, the results of the method of section IV are presented. The finite-element mesh for the domain of Fig. 1 containing 16-fourth-order curved quadrilateral elements is shown in Fig. 2. All the calculations of this section have been performed in a laptop computer equipped with a Core2Due CPU of 2.5GHz clock and 2-GBs of RAM under a 32 bit operating system.

The convergence pattern of the method of section IV is shown in Fig. 3, which is completely the same for both the product to sum approach and the direct numerical integration approach. This shows that the product to sum approach in matrix filling does not have any effect on the accuracy and the convergence of the FEM solver. The computational times for both methods are listed in Table I. As can be seen from the table, the matrix filling time in the product to sum approach is relatively small compared with the direct numerical integration approach especially when the order of the bases increases. It is also clear from the table that the matrix filling time is a small portion of the total computational time in the product to sum approach, which is in contrast with the direct numerical integration approach.

## VI. CONCLUSION

A product to sum approach has been proposed for a fast and efficient matrix filling in a high order hierarchical curved inhomogeneous FEM. In the conversion of a product of polynomials to a sum, the ideal choice is the ordinary monomials. However, the ordinary monomials suffer from the ill-conditioned finite-element matrices. The Chebyshev polynomials stand at the second rank. These polynomials were used to check the method in solving a 2-D curl-curl equation and the improvement in the efficiency of the method was verified.